\documentclass[10pt]{article}

\usepackage{graphicx} 
\usepackage{subcaption}
\usepackage[utf8]{inputenc}
\usepackage[
    a4paper,
    textwidth=5.125in,
    textheight=8.25in,
    centering
]{geometry}

\usepackage[centertags]{amsmath}
\usepackage{amssymb,amsthm}
\usepackage{amsfonts}
\usepackage[colorlinks=true,linkcolor=blue,citecolor=blue]{hyperref}
\usepackage{url}
\usepackage{mdframed}
\usepackage{multirow}
\usepackage{enumitem}
\usepackage{tikz}

\usepackage{booktabs}
\usepackage{array}
\usepackage{algorithm}
\usepackage{algpseudocode}

\newtheorem{theorem}{Theorem}
\newtheorem{lemma}[theorem]{Lemma}
\newtheorem{proposition}[theorem]{Proposition}

\newtheorem{corollary}[theorem]{Corollary}

\theoremstyle{definition}
\newtheorem{definition}[theorem]{Definition}

\newtheorem{example}[theorem]{Example}
\newtheorem{remark}[theorem]{Remark}

\title{Finite Recovery by Extrapolation for the Douglas-Rachford Method on Affine Subspaces}
\author{Jordan Collard \and Scott B. Lindstrom}

\newcommand{\dia}{\mathrm{diag}}

\newcommand{\R}{\mathbb{R}}
\newcommand{\spa}{{\operatorname {span}}}
\newcommand{\Fix}{{\rm Fix}}

\newcommand{\Rot}{{\rm Rot}}

\newcommand{\rank}{{\rm rank}}
\newcommand{\gmres}{\textnormal{DR-GMRES}}
\newcommand{\argmin}{{\operatorname {argmin}}}
\newcommand{\range}{{\operatorname {ran}}}
\newcommand{\LeSE}{{\mathrm{LeSE}}}
\newcommand{\aff}{{\operatorname{aff}}}
\newcommand{\LSM}{\mathrm{LSM}}

\newcommand{\RaDLSM}{\mathrm{RaDLSM}}
\newcommand{\RaDLeSE}{\mathrm{RaDLeSE}}

\begin{document}

\maketitle

\begin{abstract}
    We completely describe the dynamics of the Douglas--Rachford (DR) method for solving affine subspace feasibility. We characterise the affine hull of the resulting iterates, and we show how the corresponding search space's dimension is an explicit function of the number of distinct nontrivial principal angles. This dimension can be much smaller than the ambient dimension, helping to explain why extrapolation can be so effective.

    Using this characterisation, we extend the Lyapunov surrogate method (LSM) to higher extrapolation orders and prove that it recovers the DR limit exactly from finitely many iterates. We show that LSM and a least-squares extrapolation (LeSE) share the affine search spaces used by the generalized minimal residual method (GMRES) applied to the DR fixed-point equation (DR-GMRES). The dimension characterisation gives bounds on the number of steps required for exact recovery by all three methods. We develop adaptive implementations that estimate the search-space dimension from differences between successive DR iterates and use this estimate to select the extrapolation order. Numerical experiments show how ill-conditioning affects recovery accuracy, and we identify cases where adaptive extrapolation runs faster than DR-GMRES.
\end{abstract}

\section{Introduction}
The two-set feasibility problem is to find $x\in A\cap B$, where $A,B\subseteq\R^n$ are constraint sets. Such problems arise in image recovery, phase retrieval, graph colouring, among many other contexts \cite{aragon2018solving,bauschke2002phase,combettes1996convex,lindstrom2021survey}. Splitting projection methods approach these problems through the individual projections onto the constraint sets, which may be easier to compute than a projection onto their intersection. Their theory and applications in convex feasibility were reviewed in \cite{bauschke1996projection}.

Among these methods, the Douglas-Rachford method (DR) has received substantial attention in both convex and nonconvex settings. Originally introduced for the numerical solution of heat-conduction problems \cite{douglas1956numerical}, DR was developed as an operator splitting method for monotone inclusions by Lions and Mercier \cite{lions1979splitting}. For feasibility problems involving nonempty closed convex sets $A,B\subseteq\R^n$, the DR iteration from an initial point $y_0\in\R^n$ is defined by
\[
    y_{k+1}:=T_{A,B}y_k \qquad T_{A,B}:=\frac12(I+R_BR_A),
\]
where $P_C$ denotes the nearest-point projection onto $C$ and $R_C:=2P_C-I$ is the corresponding reflection. If $A$ and $B$ are convex with $A\cap B\neq \varnothing$, then $T_{A,B}$ is firmly nonexpansive, and $(y_k)$ converges to a fixed point $\bar y\in \Fix(T_{A,B})$, whose shadow $P_A\bar y\in A\cap B$ solves the feasibility problem. DR has also been applied successfully to nonconvex feasibility problems, although convergence is no longer guaranteed by the convex theory. The development of DR, its convergence theory, and its applications in convex and nonconvex settings were surveyed in \cite{lindstrom2021survey}.

Affine subspaces are important both as constraint sets and as local models for more complicated feasibility problems. Systems of linear equality constraints define affine subspaces, and their intersection represents the points satisfying the combined constraints. The associated best-approximation problem seeks the nearest such point to a prescribed initial point. Affine subspace analysis also informs local convergence theory more broadly: local reductions to affine subspace feasibility were used in the analysis of sparse affine feasibility problems \cite{hesse2014sparse}.

In a Hilbert-space setting, Bauschke et al.\ \cite{BAUSCHKE201463} established the best-approximation property for two intersecting affine subspaces:
\[
    y_k\longrightarrow P_{\Fix(T_{A,B})}y_0,\qquad P_Ay_k\longrightarrow P_{A\cap B}y_0,
\]
and identified the sharp convergence rate as the cosine of the Friedrichs angle. In finite dimensions, the principal angles between the direction subspaces \cite{bjorck1973numerical,deutsch2001best} give an explicit geometric description: on each nontrivial principal plane, the DR operator acts as a rotation combined with a contraction, causing iterates to spiral. The analysis was extended to affine subspaces that need not intersect, including the limiting behaviour of the shadow sequences \cite{bauschke2016douglas}.

This geometric description provides a basis for developing and analysing modifications of DR. In particular, a small Friedrichs angle gives a convergence factor close to one, motivating acceleration methods for DR. Bauschke et al. \ \cite{bauschke2016optimal} analysed relaxation-based acceleration by determining optimal convergence rates for generalised DR methods on two subspaces. Geometric constructions and extrapolation take a different approach, using information from several computed points to construct an improved approximation.

Behling, Bello-Cruz, and Santos introduced the circumcentered-reflection method (CRM), which replaces the averaging step of DR by the circumcenter of $y$, $R_Ay$, and $R_BR_Ay$ \cite{behling2018circumcentering}. For intersecting affine subspaces $A$ and $B$, the CRM point associated with $y$ is the orthogonal projection of any $z\in A\cap B$ onto $\operatorname{aff}\{y,R_Ay,R_BR_Ay\}$. CRM was subsequently extended to finitely many affine subspaces and to block-wise variants \cite{behling2018linear, behling2020block}. More recently, a sharp contraction bound was established for two-subspace CRM initialised in $B$ \cite{BelloCruz2026CRMrate}. This bound depends on both the Friedrichs angle and the largest principal angle, whereas the sharp convergence rate for DR depends only on the Friedrichs angle.

Lyapunov functions were used to establish convergence of DR in certain nonconvex settings \cite{benoist2015douglas, dao2019lyapunov}. Their geometry was subsequently used to relate circumcentering to spiralling DR iterates and to develop a framework for extrapolation based on minimising Lyapunov surrogates \cite{lindstrom2022computable}. Within this framework, the Lyapunov surrogate method (LSM) was introduced to construct an extrapolated point from successive iterates. If the underlying iteration remains in an affine subspace, the LSM extrapolates also belong to that subspace.

For two lines, DR exhibits a single rotation-contraction component, and LSM can recover a fixed point exactly. An extension using more iterates was also suggested for more general affine feasibility problems \cite{lindstrom2022computable}. For more general affine subspaces, several such components may be present. This raises the question of whether the extended method can still recover the DR limit exactly and, if so, how the number of required iterates depends on the geometry of the problem.

The use of successive iterates to approximate the limit of a sequence is also central to classical vector extrapolation. Cabay and Jackson introduced minimal polynomial extrapolation (MPE), which exploits polynomial relations among successive differences in affine linear sequences \cite{cabay1976polynomial}. A systematic treatment of MPE and related extrapolation methods was provided in \cite{sidi2017vector}. Poon and Liang connected vector extrapolation with the geometry of operator splitting methods \cite{poon2019trajectory,poon2020geometry}. Their approach fits a linear recurrence to observed differences and uses it to predict subsequent displacements. We refer to the resulting infinite-step extrapolation, which underlies A$^3$DMM, as least-squares extrapolation (LeSE). This construction is equivalent to MPE through a shift in the iterates used.

Krylov subspace methods provide another way to use the information generated by linear iteration. The fixed-point equation can be written as a linear system, to which the generalized minimal residual method (GMRES) \cite{saad1986gmres} can be applied. We refer to GMRES applied to the DR fixed-point equation as DR-GMRES. For equality-constrained strongly convex quadratic problems, the affine structure of the ADMM iteration was exploited to accelerate convergence using GMRES \cite{zhang2018gmres}. An Anderson extrapolation method was developed for DR applied to convex optimisation problems with linear equality constraints \cite{fu2020anderson}. It determines extrapolation coefficients from past fixed-point residuals and uses stabilisation to obtain global convergence guarantees. Walker and Ni established an essential equivalence between untruncated Anderson acceleration and GMRES for linear fixed-point problems under suitable conditions, connecting extrapolation based on previous iterates with residual minimisation over Krylov subspaces \cite{walker2011anderson}. These connections motivate examining whether the different methods search the same affine spaces and what this implies for exact recovery.

\subsection{Contributions}

For two intersecting affine subspaces, we use the principal-angle decomposition to characterise the affine space explored by DR and the dynamics within that space (Theorem~\ref{thm:dim_reduction}). Its dimension depends on the distinct principal angles between the subspaces and the initial point, and can be much smaller than the ambient dimension.

This characterisation allows us to establish finite recovery by an extension of LSM. We show that the LSM extrapolate is the closest point to the DR limit within the affine hull of the iterates used in its construction (Lemma~\ref{lem:LT_projection}). Once sufficiently many iterates have been generated to span the full search space, the extrapolate therefore recovers the limit exactly (Theorem~\ref{thm:finite_LT}). The geometry determines how many iterates suffice for this recovery.

The affine spaces generated by successive DR iterates also coincide with the affine Krylov search spaces used by DR-GMRES (Proposition~\ref{prop:common_search_space}). LSM and the LeSE select approximations from these same spaces using different criteria. This common structure gives finite recovery bounds for all three methods and shows that their full-order constructions recover the same DR limit.

To implement extrapolation without prior knowledge of the required dimension, we develop rank-adaptive methods that select their order from observed DR differences. We also introduce a QR formulation to improve the conditioning of the LSM solve. Numerical experiments examine the limitations imposed by near-linear dependence among the differences. In the tested problems, DR-GMRES requires the fewest DR operator evaluations, but the smaller systems used by the adaptive methods can yield shorter runtimes when search-space dimensions are larger.

\noindent \textbf{Outline.} After introducing the notation and preliminaries in Section~\ref{sec:preliminaries}, we characterise the DR difference dynamics and determine the search-space dimension in Section~\ref{sec:dr_dynamics_affine}. We then present LSM extrapolation, prove its exact recovery property, and describe its rank-adaptive implementation in Section~\ref{sec:LT}. In Section~\ref{sec:Krylov}, we establish the connections with LeSE and DR-GMRES. We present the numerical experiments in Section~\ref{sec:numerical_experiments} and discuss the findings and future directions in Section~\ref{sec:conclusion}.

\color{black}
\section{Preliminaries}\label{sec:preliminaries}
Throughout, \(\mathbb R^n\) is equipped with the standard inner product
\(\langle\cdot,\cdot\rangle\) and induced norm
\(\|x\|=\sqrt{\langle x,x\rangle}\). We write
\[
\delta_{ij}:=
\begin{cases}
1, & i=j,\\
0, & i\neq j,
\end{cases}
\]
for the Kronecker delta. 

Let $A=a+U$ and $B=b+V$ be affine subspaces of $\mathbb R^n$, where $a,b\in \mathbb R^n$ and $U,V\subseteq\mathbb R^n$ are linear subspaces. Throughout, we assume that $U\ne\{0\}, V\ne \{0\},$ and $A\cap B\ne \varnothing$. For a nonempty closed convex set $C\subseteq\mathbb R^n$, let $P_C$ denote the metric projection onto $C$, and define the reflection $R_C:=2P_C-I$.

Let $L:\mathbb R^n\to\mathbb R^n$ be linear. A subspace $E\subseteq\mathbb R^n$ is said to be $L$-invariant if $L(E)\subseteq E$. In this case, let $L|_E:E\to E$ denote the restriction of $L$ to $E$, and let $\mu_{L|_E}$ denote its minimal polynomial. We shall use the standard facts that if a polynomial $q$ satisfies $q(L|_E)=0$, then $\mu_{L|_E}$ divides $q$, and the roots of $\mu_{L|_E}$ are precisely the eigenvalues of $L|_E$. When $E=\{0\}$, we adopt the convention that $\mu_{L|_E}=1$, and hence $\deg\mu_{L|_E}=0$.

The notion of principal angles between subspaces dates back to Jordan \cite{jordan1875essai}. We use the recursive finite-dimensional formulation of Björck and Golub \cite{bjorck1973numerical}. Principal angles encode the geometry of two subspaces and determine both the classical convergence theory of DR and the principal-plane decomposition used below.

\begin{definition}[Principal angles and vectors]\label{def:principal_angles}
    Let $U,V\subseteq\mathbb R^n$ be nonzero linear subspaces and set $
        q:=\min\{\dim U,\dim V\}$.
    
    For $k=1,\dots,q$, having chosen $u_i$ and $v_i$ for $i<k$, choose
    unit vectors $u_k\in U$ and $v_k\in V$ attaining
    \[
        \max_{u\in U}\max_{v\in V}\langle u,v\rangle
    \]
    subject to
    \[
        \|u\|=\|v\|=1,
        \qquad 
        \langle u,u_i\rangle
        =
        \langle v,v_i\rangle
        =0
        \qquad (i<k).
    \]
    The $k$th \emph{principal angle} $\theta_k\in[0,\tfrac{\pi}{2}]$ is defined by
    \[
        \cos\theta_k:=\langle u_k,v_k\rangle,
    \]
    and $u_k$ and $v_k$ are called corresponding \emph{principal vectors}. The principal angles satisfy $0\leq\theta_1\leq\cdots\leq\theta_q\leq\frac{\pi}{2}.
    $
\end{definition}

The following standard result makes this geometric information explicit through mutually orthogonal two-dimensional principal planes; see \cite{halmos1969two,wedin1983angles,tian2019relationships}.

\begin{lemma}[Principal planes]\label{lem:princ_planes}
    Let $\theta_1,\dots,\theta_q$ and corresponding principal vectors
    $u_i\in U$ and $v_i\in V$ be as in Definition~\ref{def:principal_angles}.
    Then
    \[
        \langle u_i,v_j\rangle=\delta_{ij}\cos\theta_i.
    \]
    For each $i$ such that $0<\theta_i<\tfrac{\pi}{2}$, define
    \[
        \widetilde v_i
        :=\frac{v_i-\cos\theta_i\,u_i}{\sin\theta_i},
        \qquad
        \Pi_i:=\spa\{u_i,\widetilde v_i\}.
    \]
    Then $\widetilde v_i\in U^\perp$, and
    $(u_i,\widetilde v_i)$ is an orthonormal basis of $\Pi_i$.
    The subspace $\Pi_i$ is called the corresponding \emph{principal plane}.
\end{lemma}
The associated boundary cases and orthogonality relations yield the canonical decomposition of the ambient space; see \cite[Section 2 and Appendix 2]{wedin1983angles} and \cite[Section 2]{tian2019relationships}.

\begin{lemma}[Principal-angle structure]
\label{lem:principal_structure}
    The principal vectors and principal planes satisfy the following.
    \begin{enumerate}
        \item \emph{Boundary cases.}
        \begin{itemize}
            \item If $\theta_i=0$, then
            $u_i=v_i\in U\cap V$.
            \item If $\theta_i=\tfrac{\pi}{2}$, then
            $u_i\in U\cap V^\perp$ and
            $v_i\in U^\perp\cap V$.
        \end{itemize}

        \item \emph{Orthogonality.}
        The principal planes $\Pi_i$ are mutually orthogonal and are
        orthogonal to
        \[
            U\cap V,\qquad U\cap V^\perp,\qquad
            U^\perp\cap V,\qquad U^\perp\cap V^\perp.
        \]

        \item \emph{Orthogonal decompositions.}
        \begin{align*}
            U
            &=(U\cap V)
              \oplus
              \bigoplus_{0<\theta_i<\tfrac{\pi}{2}}
              \spa\{u_i\}
              \oplus
              (U\cap V^\perp),\\
            V
            &=(U\cap V)
              \oplus
              \bigoplus_{0<\theta_i<\tfrac{\pi}{2}}
              \spa\{v_i\}
              \oplus
              (U^\perp\cap V),\\
            \mathbb R^n
            &=(U\cap V)
              \oplus(U\cap V^\perp)
              \oplus(U^\perp\cap V)
              \oplus(U^\perp\cap V^\perp)
              \oplus
              \bigoplus_{0<\theta_i<\tfrac{\pi}{2}}\Pi_i.
        \end{align*}
    \end{enumerate}
\end{lemma}

We shall use the following standard direct-sum property of minimal polynomials.
\begin{lemma}[Minimal polynomial on invariant direct sums]\label{lem:lcm}
    Let $X$ be a finite-dimensional vector space over $\mathbb{F}\in \{\mathbb{R},\mathbb{C}\}$ and let $T\in \mathcal{L}(X)$. Suppose that $X=X_1\oplus\dots\oplus X_m$ where each $X_i$ is $T$-invariant. For each $i$ let $T_i:=T|_{X_i}\in\mathcal{L}(X_i)$
    be the restriction of $T$ to $X_i$. Then $\mu_T=\mathrm{lcm}\big(\mu_{T_1},\dots,\mu_{T_m}\big)$, where $\mathrm {lcm}$ denotes the least common multiple of the polynomials $\mu_{T_1},\dots,\mu_{T_m}$.
    
\end{lemma}

\begin{proof}
    Let $\mu_i:=\mu_{T_i}$ for $i=1,\dots,m$. Fix a polynomial $p\in\mathbb{F}[x]$. Since each $X_i$ is $T$-invariant, we have $p(T)(X_i)\subseteq X_i$ and $p(T)|_{X_i}=p(T_i)$ for $i=1,\dots,m$. Since $X=X_1\oplus\dots\oplus X_m$, the operator $p(T)$ is the zero operator on $X$ if and only if its restriction to each $X_i$ is zero. Hence we obtain the following
    \[
    p(T)=0\iff \forall i, \,p(T_i)=0\iff \forall i,\, \mu_i\mid p.
    \]
    
    By definition, the minimal polynomial $\mu_T$ is the unique monic polynomial of smallest degree satisfying $p(T)=0$. Thus $\mu_T$ must be the smallest monic polynomial divisible by all the $\mu_i$, which is exactly their least common multiple. Thus $\mu_T=\mathrm{lcm}(\mu_1,\dots,\mu_m)$.
\end{proof}

\section{Douglas--Rachford difference dynamics on affine subspaces}\label{sec:dr_dynamics_affine}
We will first recall the Douglas--Rachford operator and show that the difference sequence of an affine-subspace DR is generated by repeated application of the linear operator associated with the corresponding direction subspaces. We then describe the action of this linear operator on the principal-angle decomposition of Lemma~\ref{lem:principal_structure}. This blockwise representation provides the basis for the minimal polynomial analysis and difference-space characterisation developed below.

\begin{definition}[Douglas--Rachford operator]
    Let \(C,D\subseteq \mathbb R^n\) be nonempty closed convex sets. The Douglas--Rachford operator associated with $C,D$ is
    \[
    T_{C,D}:=\frac12(I+R_DR_C).
    \]
\end{definition}

Although our main results concern affine subspaces, the difference dynamics of the affine DR iteration are governed by the corresponding direction subspaces. Indeed, for two affine subspaces $A=a+U$, $B=b+V$, the operator $T_{A,B}$ has the form $T_{A,B}x=Lx+c$, where $L:=T_{U,V}$, for some $c\in \R^n$. Consequently, for any DR sequence $y_{k+1}:=T_{A,B}y_k$, the difference vectors $d_k:=y_k-y_{k-1}$ satisfy
\begin{equation}\label{eq:difference_recurrence}
    d_{k+1}=T_{A,B}y_k-T_{A,B}y_{k-1}=L(y_k-y_{k-1})=Ld_k.
\end{equation}
Hence $d_k=L^{k-1}d_1$ for every $k\ge1$. Thus, the principal-angle structure of the linear operator $L=T_{U,V}$ determines the invariant subspace generated by an affine DR difference sequence. We first analyse the blockwise action and minimal polynomial of $L$, and then use this structure to characterise the space spanned by affine DR differences.

The next result gives the blockwise action of \(L:=T_{U,V}\) on the components of the principal-angle decomposition in Lemma~\ref{lem:principal_structure}. It follows from the standard principal-plane representation of two subspaces; see, e.g., \cite{BAUSCHKE201463,bauschke2016douglas}. For completeness, a proof is included in Appendix~\ref{app:dr_proof}.
\begin{proposition}[Blockwise action of DR]\label{prop:angles}
    Let $\Pi_i$ be a principal plane corresponding to a principal angle $\theta_i\in (0,\pi/2)$. Then $\Pi_i$ is $L$-invariant and, with respect to the orthonormal basis $(u_i,\tilde v_i)$, 
    \[
        L|_{\Pi_i}=\cos\theta_i\Rot(\theta_i),\qquad \Rot(\theta_i)=\begin{bmatrix}
            \cos\theta_i & -\sin\theta_i\\
            \sin\theta_i & \cos\theta_i
        \end{bmatrix}.
    \]
    Moreover, 
    \[
    L|_{U\cap V}=I,\quad L|_{U^\perp\cap V^\perp}=I,\quad L|_{U\cap V^\perp}=0,\quad L|_{U^\perp \cap V}=0.
    \]
\end{proposition}

Motivated by Proposition~\ref{prop:angles}, set
\[
    \mathcal W:=\bigoplus_{0<\theta_i<\frac{\pi}{2}}\Pi_i,
    \qquad
    N:=(U\cap V^\perp)\oplus(U^\perp\cap V),
\]
and define the trivial and nontrivial components 
$$
F:=(U\cap V)\oplus(U^\perp\cap V^\perp),\qquad G:=N\oplus \mathcal W,
$$ 
respectively. Then $\mathbb R^n=F\oplus G$, $ L|_F=I$, and $L|_N=0$. Since each $\Pi_i$ is $L$-invariant, so are $\mathcal W$ and $G=N\oplus \mathcal W$. Moreover, $L|_{\Pi_i}$ is invertible and has no nonzero fixed points for every $i$. Consequently, $\Fix(L)=F$ and $\ker L=N$.

The following convergence result is standard; see \cite[Corollary~4.5]{BAUSCHKE201463}. We include a short proof adapted to the principal-angle decomposition used here.

\begin{proposition}[Limit of affine DR iterates]\label{prop:affine_DR_limit}
For any $y_0\in \R^n$, let $y_{k+1}=T_{A,B}y_k$ and set $\bar y:=P_{\Fix(T_{A,B})}y_0$. Then
\[
y_k-\bar y=L^k(y_0-\bar y),\qquad y_0-\bar y\in G,
\]
and consequently $y_k\longrightarrow P_{\Fix(T_{A,B})}y_0$.
\end{proposition}

\begin{proof}
    Since $\bar y\in\Fix(T_{A,B})$, for every $x\in\R^n$ we have $T_{A,B}x-T_{A,B}\bar y=L(x-\bar y)$. It follows inductively that $y_k-\bar y=L^k(y_0-\bar y)$. Moreover, 
    \[
    \begin{aligned}
        x\in \Fix(T_{A,B})
        \iff L(x-\bar y)=x-\bar y
        \iff x-\bar y\in F.
    \end{aligned}
    \]
    Hence $\Fix(T_{A,B})=\bar y+F.$ Because $\bar y$ is the projection of $y_0$ onto this affine subspace, $y_0-\bar y\in F^\perp =G$.

    On $N$, the operator $L$ is zero, while on each principal plane $\Pi_i$ associated with $\theta_i\in (0,\pi/2)$,
    \[
        L|_{\Pi_i}=\cos(\theta_i)\operatorname{Rot}(\theta_i).
    \]
    Consequently, $\|L^k|_{\Pi_i}\|=\cos^k(\theta_i)\longrightarrow 0.$ Therefore
    \[
        y_k-\bar y=L^k(y_0-\bar y)\longrightarrow 0,
    \]
    which proves the result.
\end{proof}

Let $\Theta:=\{\theta_i:0<\theta_i<\pi/2\}$ be the set of distinct principal angles lying in $(0,\pi/2)$, and let $m:=|\Theta|$. For each $\theta\in\Theta$, define $\mathcal W_\theta:=\bigoplus_{\{i:\theta_i=\theta\}}\Pi_i$. Then $\mathcal W=\bigoplus_{\theta\in \Theta}\mathcal W_\theta$.
\begin{proposition}[Minimal-polynomial degree on $G$]
\label{prop:deg_T}
    Let $\ell_G:=\deg\mu_{L|_G}$. Then
    \[
        \ell_G=
        \begin{cases}
            2m,   &\text{if }N=\{0\},\\
            2m+1, &\text{otherwise}.
        \end{cases}
    \]
\end{proposition}

\begin{proof}
For each $\theta\in\Theta$ and each principal plane $\Pi_i$ with $\theta_i=\theta$, Proposition~\ref{prop:angles} gives $L|_{\Pi_i}=\cos\theta\,\Rot(\theta)$. The complex eigenvalues of this restriction are $\lambda_\theta^\pm=\cos\theta\,e^{\pm\mathrm{i}\theta}$. Since $0<\theta<\pi/2$, these eigenvalues are distinct nonreal complex conjugates, so the minimal polynomial of $L|_{\Pi_i}$ must have both as roots. Hence, whenever $\theta_i=\theta$,
\[
\mu_{L|_{\Pi_i}}(x)
=(x-\lambda_\theta^+)(x-\lambda_\theta^-)
=x^2-2\cos^2(\theta)x+\cos^2(\theta)
=:p_\theta(x).
\]

Since $\mathcal W_\theta=\bigoplus_{\{i:\theta_i=\theta\}}\Pi_i$, Lemma~\ref{lem:lcm} applied to $L|_{\mathcal W_\theta}$ gives
\[
\mu_{L|_{\mathcal W_\theta}}
=\operatorname{lcm}\{\mu_{L|_{\Pi_i}}:\theta_i=\theta\}
=p_\theta.
\]
If $\theta,\varphi\in\Theta$ are distinct, then $p_\theta$ and $p_\varphi$ have no common complex root. Indeed, the roots of $p_\theta$ have modulus $\cos\theta$, while the roots of $p_\varphi$ have modulus $\cos\varphi$, and $\cos\theta\neq\cos\varphi$. Thus $p_\theta$ and $p_\varphi$ are coprime. Applying Lemma~\ref{lem:lcm} to $\mathcal W=\bigoplus_{\theta\in\Theta}\mathcal W_\theta$ therefore gives
\[
\mu_{L|_\mathcal W}(x)=\prod_{\theta\in\Theta}p_\theta(x),
\]
where this identity also holds when $\Theta=\varnothing$, by the convention $\mu_{L|_{\{0\}}}=1$. Hence $\deg\mu_{L|_\mathcal W}=2m$.

If $N=\{0\}$, then $G=\mathcal W$, so $\ell_G=2m$. Suppose now that $N\neq\{0\}$. Since $G=N\oplus \mathcal W$ and $L|_N=0$, we have $\mu_{L|_N}(x)=x$. Lemma~\ref{lem:lcm} gives $\mu_{L|_G}(x)=\operatorname{lcm}\bigl(x,\mu_{L|_\mathcal W}(x)\bigr).$ Since $p_\theta(0)=\cos^2\theta\neq0$ for every $\theta\in\Theta$, the polynomial $x$ is coprime with $\mu_{L|_\mathcal W}$. Therefore $\mu_{L|_G}(x)=x\mu_{L|_\mathcal W}(x)$, and consequently $\ell_G=2m+1$.
\end{proof}

The preceding proposition describes the degree of the minimal polynomial on the full nontrivial component $G$. However, the initial DR difference need not have a nonzero component in every invariant block of $G$. Before determining the dimension of the subspace generated by the difference sequence, we record several basic properties of affine DR difference sequences.

\begin{lemma}[Properties of the affine DR difference space]
\label{lem:difference_space_properties} 
Let $\mathcal S:=\operatorname{span}\{d_1,d_2,\ldots\},$ and $ s:=\dim\mathcal S.$ 
Then the following hold:
\begin{enumerate}
    \item[(i)] $d_k=L^{k-1}d_1$ for every $k\geq1$;
    \item[(ii)] $\mathcal S\subseteq G$ and $\mathcal S$ is $L$-invariant;
    \item[(iii)] $y_k\in y_0+\mathcal S$ for every $k\geq0$;
    \item[(iv)] if $s\geq1$, then $d_1,\ldots,d_s$ form a basis of $\mathcal S$.
\end{enumerate}
\end{lemma}

\begin{proof}
For \textup{(i)}, the affine difference recurrence \eqref{eq:difference_recurrence} gives $d_{k+1}=Ld_k$. It follows inductively that $d_k=L^{k-1}d_1$ for every $k\geq1$.
 
For \textup{(ii)}, choose $x^\star\in A\cap B$. Since $T_{A,B}x^\star=x^\star$, we have $d_1=T_{A,B}y_0-y_0=(L-I)(y_0-x^\star).$ Since $L|_F=I$ and $G=F^\perp$ is $L$-invariant, $d_1=(L-I)P_G(y_0-x^\star)\in G.$ It follows from $d_k=L^{k-1}d_1$ that $d_k\in G$ for every $k\geq1$, and therefore $\mathcal S\subseteq G$. Moreover, $Ld_k=d_{k+1}\in\mathcal S$, so $\mathcal S$ is $L$-invariant.

For \textup{(iii)}, the telescoping identity $y_k-y_0=\sum_{j=1}^k d_j$ shows that $y_k-y_0\in\mathcal S$, and hence $y_k\in y_0+\mathcal S$ for every $k\geq0$.

For \textup{(iv)}, suppose that $s\ge 1$. Then $d_1\ne 0$, since $d_1=0$ would imply $d_k=L^{k-1}d_1=0$ for every $k\ge 1$, and hence $s=0$. Suppose, to the contrary, that $d_1,\dots,d_s$ are linearly dependent. Since $d_1\ne 0$, there exists a smallest $j\in \{2,\dots,s\}$ such that $d_j\in\mathcal K:=\spa\{d_1,\dots,d_{j-1}\}$. The recurrence $Ld_i=d_{i+1}$ gives $Ld_i\in\mathcal K$ for $i=1,\dots,j-2$, while $Ld_{j-1}=d_j\in\mathcal K$. Hence $L\mathcal K\subseteq\mathcal K$. It follows inductively that $d_k\in \mathcal K$ for every $k\ge j$. Therefore $\mathcal S=\mathcal K=\spa\{d_1,\dots,d_{j-1}\}$, which yields $s=\dim\mathcal S\le j-1<s$, a contradiction. Thus $d_1,\dots,d_s$ are linearly independent. Since $\dim\mathcal S=s$, they form a basis of $\mathcal S$.
\end{proof}

The recurrence $d_k=L^{k-1}d_1$ shows that $S$ is the cyclic subspace generated by $d_1$ under $L$. For $s\ge 1$, Lemma~\ref{lem:difference_space_properties} gives $S=\spa\{d_1,Ld_1,\dots,L^{s-1}d_1\}$. A space spanned by successive applications of a linear operator to a vector is called a Krylov subspace. We discuss this connection further in Section~\ref{subsec:GMRES}. We now identify the invariant blocks activated by the initial difference and determine the exact dimension of the resulting difference space.
\begin{theorem}[Dimension of the affine DR difference space]
\label{thm:dim_reduction}
Let $\mathcal S$ and $s$ be as in Lemma~\ref{lem:difference_space_properties}. For each $\theta\in\Theta$, set $d_{1,\theta}:=P_{W_\theta}d_1,\,d_{1,N}:=P_Nd_1,$ and define $\Theta_{\rm act}:=\{\theta\in\Theta:d_{1,\theta}\neq\mathbf{0}\}.$ Let $\delta_N:=1$ if $d_{1,N}\neq\mathbf{0}$, and let $\delta_N:=0$ otherwise. Then
\[
\mathcal S
=
\operatorname{span}\{d_{1,N}\}
\oplus
\bigoplus_{\theta\in\Theta_{\rm act}}
\operatorname{span}\{d_{1,\theta},Ld_{1,\theta}\}.
\]
Consequently, $s=2|\Theta_{\rm act}|+\delta_N\leq\ell_G$.
\end{theorem}

\begin{proof}
By Lemma~\ref{lem:difference_space_properties}, we have $d_1\in G$ and $\mathcal S=\operatorname{span}\{d_1,Ld_1,L^2d_1,\ldots\}.$ If $d_1=\mathbf{0}$, then $\mathcal S=\{0\}$, $\Theta_{\rm act}=\varnothing$, and $\delta_N=0$, so the result is immediate. Hence suppose that $d_1\neq\mathbf{0}$.

The orthogonal decomposition $G=N\oplus\bigoplus_{\theta\in\Theta}W_\theta$ gives
\[
d_1=d_{1,N}+\sum_{\theta\in\Theta}d_{1,\theta}.
\]
By Proposition~\ref{prop:angles}, $L|_N=0$. For each $\theta\in\Theta$, the polynomial $p_\theta(t)=t^2-2\cos^2(\theta)t+\cos^2(\theta)$ annihilates $L|_{W_\theta}$, i.e., $L^2-2\cos^2(\theta)L+\cos^2(\theta)I=0$ on $W_\theta$. Hence $L^2d_{1,\theta}=2\cos^2(\theta)Ld_{1,\theta}-\cos^2(\theta)d_{1,\theta}.$ Thus $\operatorname{span}\{d_{1,\theta},Ld_{1,\theta}\}$ is $L$-invariant. Similarly, since $L|_N=0$, we have $Ld_{1,N}=\mathbf{0}$, so $\spa\{d_{1,N}\}$ is $L$-invariant. Therefore,
\[
\mathcal H
:=
\operatorname{span}\{d_{1,N}\}
\oplus
\bigoplus_{\theta\in\Theta_{\rm act}}
\operatorname{span}\{d_{1,\theta},Ld_{1,\theta}\}
\]
is $L$-invariant and contains $d_1$. Since $\mathcal S=\operatorname{span}\{d_1,Ld_1,L^2d_1,\ldots\},$ it follows that $\mathcal S\subseteq\mathcal H$.

For each $\theta\in\Theta_{\rm act}$, the vectors $d_{1,\theta}$ and $Ld_{1,\theta}$ are linearly independent. Otherwise, $Ld_{1,\theta}=\lambda d_{1,\theta}$ for some $\lambda\in\mathbb R$, and $p_\theta(L)d_{1,\theta}=\mathbf{0}$ would imply $p_\theta(\lambda)=0$. This is impossible because the discriminant of $p_\theta$ is $-4\cos^2(\theta)\sin^2(\theta)<0.$ Since the decomposition of $G$ is orthogonal, it follows that $\dim\mathcal H=2|\Theta_{\rm act}|+\delta_N.$

It remains to show that $\mathcal S=\mathcal H$. Set
\[
r(t):=t^{\delta_N}\prod_{\theta\in\Theta_{\rm act}}p_\theta(t).
\]
On $N$, $Ld_{1,N}=\mathbf{0}$ whenever $\delta_N=1$, while on each active space $W_\theta$, $p_\theta(L)d_{1,\theta}=\mathbf{0}$. Since these spaces form an $L$-invariant direct sum, it follows that $r(L)d_1=\mathbf{0}.$

We show that $r$ divides every polynomial $q$ satisfying $q(L)d_1=\mathbf{0}$. Let $q$ be such a polynomial. Since $N$ and the spaces $W_\theta$ form an $L$-invariant direct-sum decomposition of $G$, the operator $q(L)$ acts separately on each component. Therefore, for each $\theta\in\Theta_{\rm act}$, $\mathbf{0}=P_{W_\theta}q(L)d_1=q(L)d_{1,\theta}.$ Divide $q$ by $p_\theta$ to obtain 
$$
q=a_\theta p_\theta+b_\theta,\qquad\deg b_\theta\leq1.
$$ 
Since $p_\theta(L)d_{1,\theta}=\mathbf{0}$, we have $b_\theta(L)d_{1,\theta}=\mathbf{0}$. The linear independence of $d_{1,\theta}$ and $Ld_{1,\theta}$ implies that $b_\theta=0$, and hence $p_\theta\mid q$. If $d_{1,N}\neq\mathbf{0}$, projecting onto $N$ gives $\mathbf{0}=q(L)d_{1,N}=q(0)d_{1,N},$ so $q(0)=0$ and therefore $t\mid q$.

By the proof of Proposition~\ref{prop:deg_T}, the polynomials $p_\theta$ corresponding to distinct values of $\theta$ are pairwise coprime. Moreover, $p_\theta(0)=\cos^2(\theta)\neq0$, so each $p_\theta$ is coprime with $t$. Hence the factors in $r$ are pairwise coprime, and therefore $r\mid q$. Since $r(L)d_1=\mathbf{0}$ and every polynomial $q$ satisfying $q(L)d_1=\mathbf{0}$ is divisible by $r$, the polynomial $r$ is the monic polynomial of least degree satisfying $r(L)d_1=\mathbf{0}$.

By the minimality of $r$, the vectors $d_1,Ld_1,\ldots,L^{\deg r-1}d_1$ are linearly independent. Moreover, since $r$ is monic, the identity $r(L)d_1=\mathbf{0}$ implies, by induction, that every higher power of $L$ applied to $d_1$ belongs to their span. Therefore,
\[
\dim\mathcal S
=
\deg r
=
2|\Theta_{\rm act}|+\delta_N
=
\dim\mathcal H.
\]
Since $\mathcal S\subseteq\mathcal H$, we conclude that $\mathcal S=\mathcal H$, that is,
\[
\mathcal S
=
\operatorname{span}\{d_{1,N}\}
\oplus
\bigoplus_{\theta\in\Theta_{\rm act}}
\operatorname{span}\{d_{1,\theta},Ld_{1,\theta}\}.
\]

Finally, $\Theta_{\rm act}\subseteq\Theta$, while $\delta_N=0$ if $N=\{0\}$ and $\delta_N\leq1$ otherwise. Proposition~\ref{prop:deg_T} therefore gives $s\leq\ell_G$.
\end{proof}

\begin{corollary}[Attainment of the dimension bound]
\label{cor:dimension_bound_attainment}
Under the assumptions of Theorem~\ref{thm:dim_reduction}, $s=\ell_G$ if and only if $d_{1,\theta}\neq 0$ for every $\theta\in\Theta$ and, whenever $N\neq\{0\}$, $d_{1,N}\neq 0$.
\end{corollary}

\begin{proof}
Theorem~\ref{thm:dim_reduction} gives $s=2|\Theta_{\rm act}|+\delta_N$, and $\Theta_{\rm act}\subseteq\Theta$. If $N=\{0\}$, then $\delta_N=0$ and Proposition~\ref{prop:deg_T} gives $\ell_G=2|\Theta|$. Hence $s=\ell_G$ if and only if $\Theta_{\rm act}=\Theta$. If $N\neq\{0\}$, then Proposition~\ref{prop:deg_T} gives $\ell_G=2|\Theta|+1$, while $\delta_N\in\{0,1\}$. Hence $s=\ell_G$ if and only if $\Theta_{\rm act}=\Theta$ and $\delta_N=1$. By definition, these conditions are equivalent to $d_{1,\theta}\neq0$ for every $\theta\in\Theta$ and, whenever $N\neq\{0\}$, $d_{1,N}\neq0$.
\end{proof}

Theorem~\ref{thm:dim_reduction} is the main dimension-reduction result of this section. It shows that, although DR acts on the ambient space $\mathbb R^n$, the DR iterates lie in the affine subspace $y_0+\mathcal S$. Each active principal-angle value contributes two dimensions to $\mathcal S$, while the null component contributes at most one additional dimension. Thus $s=2|\Theta_{\rm act}|+\delta_N\leq \ell_G$, and this bound is independent of the ambient dimension and of the multiplicities of the principal-angle blocks.

The following example illustrates this dimension reduction in a simple three-dimensional setting.
\begin{figure}[h]
    \centering
    \includegraphics[width=0.5\linewidth]{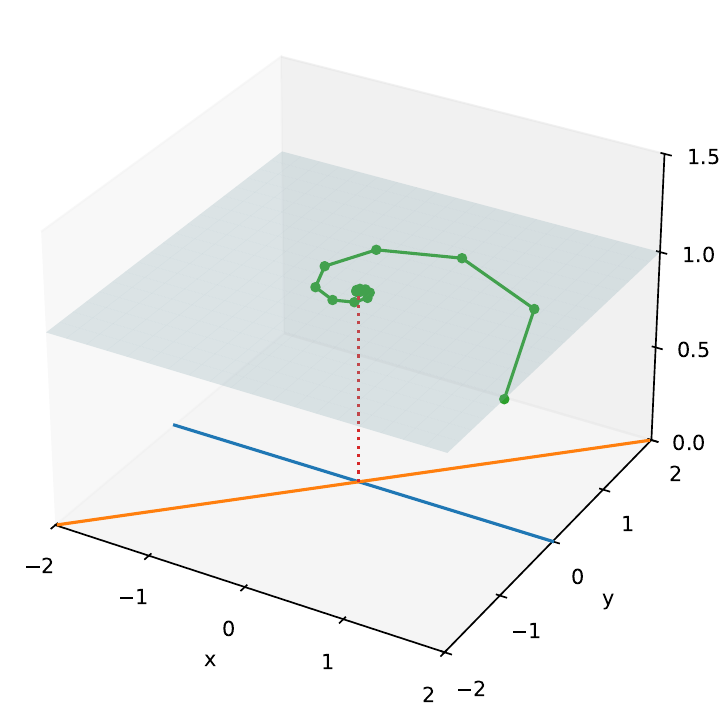}
    \caption{DR iterates for two intersecting lines in the affine plane $z=1$.}
    \label{fig:3d}
\end{figure}

\begin{example}
    Consider applying DR to two intersecting lines contained in the affine plane \(z=1\subset\mathbb R^3\), see Figure~\ref{fig:3d}. The associated direction subspaces have one nonzero principal angle $\theta\in(0,\pi/2)$ and satisfy
    \[
        U^\perp\cap V=\{0\},
        \qquad 
        U\cap V^\perp=\{0\}.
    \]
    Hence $N=\{0\}$ and Proposition~\ref{prop:deg_T} gives $\ell_G=2m=2$. For any nonstationary DR sequence, $d_1\ne0$, so the unique principal-angle block is active. Corollary~\ref{cor:dimension_bound_attainment} therefore gives $s=\ell_G=2$. Moreover, DR preserves the $z$-coordinate, so the iterates remain confined to a two-dimensional affine plane parallel to $z=1$.
\end{example}

Beyond this dimension reduction, DR iterates satisfy the following useful orthogonality relation.

\begin{proposition}[Affine iterate--difference orthogonality]
\label{prop:iterate_difference_orth}
    Let $(y_k)$ be an affine DR sequence and $d_k:=y_k-y_{k-1}$ for $k\ge 1$. Then, for every $\hat y\in \Fix(T_{A,B})$,
    \[
        \langle y_k-\hat y,d_k\rangle=0
        \qquad k\geq1.
    \]
\end{proposition}

\begin{proof}
    Fix $\hat y\in \Fix(T_{A,B})$ and define $z_k:=y_k-\hat y$ for $k\ge 0$. Since $\hat y=T_{A,B}\hat y=L\hat y+c$, we have 
    \begin{align*}
        z_k
        =
        y_k-\hat y
        =
        T_{A,B}y_{k-1}-T_{A,B}\hat y
        =
        L(y_{k-1}-\hat y)
        =
        Lz_{k-1}.
    \end{align*}
    By \cite[Proposition 3.5(iii)-(iv)]{BAUSCHKE201463}, the linear DR operator $L$ satisfies $LL^\top =L^\top L$ and $2LL^\top =L+L^\top$. Therefore
    \begin{align*}
        \|z_k\|^2
        &=
        \|Lz_{k-1}\|^2\\
        &=
        \langle L^\top Lz_{k-1},z_{k-1}\rangle\\
        &=
        \langle LL^\top z_{k-1},z_{k-1}\rangle\\
        &=
        \tfrac12\langle (L+L^\top)z_{k-1},z_{k-1}\rangle\\
        &=
        \tfrac12\big(\langle Lz_{k-1},z_{k-1}\rangle+\langle z_{k-1},Lz_{k-1}\rangle \big)\\
        &=
        \langle Lz_{k-1},z_{k-1}\rangle\\
        &=
        \langle z_k,z_{k-1}\rangle.
    \end{align*}
    It follows that 
    \[
    \langle z_k,z_k-z_{k-1}\rangle=\|z_k\|^2-\langle z_k,z_{k-1}\rangle=0.
    \]
    Finally, $z_k=y_k-\hat y$ and $z_k-z_{k-1}=y_k-y_{k-1}=d_k$. Hence $\langle y_k-\hat y,d_k\rangle =0$.
\end{proof}
Together, Theorem~\ref{thm:dim_reduction} and Proposition~\ref{prop:iterate_difference_orth} reveal both the dimension of the DR difference space and the geometry of the iterates. In the next section, we exploit these properties to construct an extrapolation scheme based on successive DR iterates with finite convergence.

\section{Lyapunov surrogate extrapolation: exact recovery and rank adaptation}\label{sec:LT}
The Lyapunov surrogate method \cite{lindstrom2022computable} constructs an extrapolated point by exploiting orthogonality relations like those in Proposition~\ref{prop:iterate_difference_orth}. We first define the $\LSM$ extrapolate and establish its projection characterisation and exact recovery property for affine DR. We then describe a rank-adaptive implementation and its QR formulation.
\subsection{LSM extrapolation and exact recovery} \label{subsec:LSM_extrap}
Given successive iterates $y_0,\dots,y_p$, we seek a point in their affine hull whose displacement from each $y_k$ is orthogonal to the corresponding differences $d_k$.
\begin{definition}[$\LSM_p$ extrapolate]\label{def:lt_p}
    Let $T:\R^n\to \R^n$, fix $y_0\in\R^n$ and define $y_{k+1}=Ty_k$. For $k\ge 1$, set $d_k:=y_k-y_{k-1}$, and $w_k:=y_k-y_0$. Let $p\ge 1$ be the extrapolation order, and suppose $d_1,\dots,d_p$ are linearly independent. Then the $\LSM_p$ extrapolate associated with $y_0,\dots,y_p$ is
    \[
    y^{\LSM_p}
    := 
    y_0+\sum_{j=1}^p\alpha_j w_j,
    \]
    where $\alpha \in \R^p$ satisfies $\langle y^{\LSM_p}-y_k,d_k \rangle =0,$ for $k=1,\dots,p$.
\end{definition}
For each $k\ge1$, define $D_k:=[d_1,\dots,d_k]$ and $W_k:=[w_1,\dots,w_k]$. For an $\LSM_p$ extrapolate, the corresponding orthogonality conditions give, for $k=1,\dots,p$,
\[
    0
    =
    \langle y_0+W_p\alpha -y_k,d_k\rangle 
    =
    d_k^\top W_p\alpha -d_k^\top w_k.
\]
Therefore, stacking these $p$ equations yields
\begin{equation}\label{eqn:LSM_coefficient_system}
    D_p^\top W_p\alpha
    =
    \begin{bmatrix}
        \langle d_1,w_1\rangle\\
        \vdots\\
        \langle d_p,w_p\rangle
    \end{bmatrix}
    =
    \dia(D_p^\top W_p),
\end{equation}
where $\dia(M)$ denotes the vector containing the diagonal entries of $M$.

The linear independence assumption in Definition~\ref{def:lt_p} ensures that the extrapolate is well defined. Indeed, since $w_k=\sum_{j=1}^kd_j$, we have $W_p=D_p\Gamma_p$, where $\Gamma_p$ is the upper-triangular matrix given by 
\[
    (\Gamma_p)_{ij}
    =
    \begin{cases}
        1, & i\le j,\\
        0, & i>j.
    \end{cases}
\]
Since $d_1,\dots,d_p$ are linearly independent, $D_p^\top D_p$ is positive definite. Moreover, $\Gamma_p$ is nonsingular, since it is upper triangular with unit diagonal. Therefore $D_p^\top W_p=(D_p^\top D_p)\Gamma_p$ is nonsingular, and hence the coefficient vector $\alpha$, and consequently the $\LSM_p$ extrapolate, are uniquely determined.

In the affine DR setting, the $\LSM$ extrapolated point admits a geometric interpretation as the orthogonal projection of the limiting fixed point $\bar y:=P_{\Fix (T_{A,B})}y_0 $ onto the affine space generated by the first $p$ difference vectors. For $p<s$, this affine space is a proper affine subspace of $y_0+\mathcal S$.

\begin{lemma}[Projection characterisation of LSM] \label{lem:LT_projection}
Let $\bar y:=P_{\Fix(T_{A,B})}y_0$. Suppose $d_1,\dots,d_p$ are linearly independent, and define $\mathcal A_p:=y_0+\spa\{d_1,\dots,d_p\}$. Then
    \[
        y^{\LSM_p}=P_{\mathcal A_p}\bar y.
    \]
\end{lemma}
\begin{proof}
    By Definition~\ref{def:lt_p}, we have $y^{\LSM_p}\in \mathcal A_p$ and 
    \[
        \langle y^{\LSM_p}-y_k,d_k\rangle =0\qquad k=1,\dots,p.
    \]
    Since $\bar y\in \Fix(T_{A,B})$, Proposition~\ref{prop:iterate_difference_orth} gives
    \[
        \langle y_k-\bar y,d_k\rangle =0\qquad k=1,\dots,p.
    \]
    Combining these identities shows
    \[
        y^{\LSM_p}-\bar y\in \spa \{d_1,\dots,d_p\}^\perp.
    \]
    Together with $y^{\LSM_p}\in \mathcal A_p$, this characterises $y^{\LSM_p}$ as $P_{\mathcal A_p}\bar y$.
\end{proof}

Lemma~\ref{lem:LT_projection} reduces the full-order analysis to determining when the limiting fixed point belongs to the affine space generated by the difference sequence.
\begin{theorem}[Exact recovery by $\LSM_s$]
\label{thm:finite_LT} 
Fix $y_0\in\mathbb R^n$. With $\mathcal S$ and $s$ as in Lemma~\ref{lem:difference_space_properties},
if $s=0$, then $y_0\in\Fix(T_{A,B})$. If $s\geq1$, let $y^{\LSM_s}$ be the $\LSM_s$ point associated with $y_0,\ldots,y_s$ as in Definition~\ref{def:lt_p}. Then $$y^{\LSM_s}=P_{\Fix(T_{A,B})}y_0.$$
\end{theorem}

\begin{proof}
    If $s=0$, then $d_1=0$, so $y_0\in \Fix(T_{A,B})$. Suppose $s\ge 1$, and set $\bar y:=P_{\Fix(T_{A,B})}y_0$. By Lemma~\ref{lem:difference_space_properties}, $d_1,\dots,d_s$ form a basis of $\mathcal S$, and $y_k\in y_0+\mathcal S$ for every $k\ge 0$. Proposition~\ref{prop:affine_DR_limit} gives $y_k\to \bar y$. Since $\mathcal S$ is a subspace of the finite-dimensional space $\R^n$, the affine space $y_0+\mathcal S$ is closed. It follows that 
    \begin{equation}\label{eq:projection_in_difference_space}
        \bar y-y_0\in \mathcal S.
    \end{equation}
    Applying Lemma~\ref{lem:LT_projection}, we obtain $y^{\LSM_s}=P_{y_0+\mathcal S}\bar y=\bar y=P_{\Fix(T_{A,B})}y_0.$
\end{proof}

Thus, full-order $\LSM$ recovers the distinguished fixed point after at most $s$ difference directions have been generated.
\begin{remark}
By \cite[Corollary~4.5]{BAUSCHKE201463}, $P_A P_{\Fix(T_{A,B})}y_0=P_{A\cap B}y_0.$ Consequently, Theorem~\ref{thm:finite_LT} shows that the shadow of the $\LSM_s$ point recovers the nearest feasible point $P_Ay^{\LSM_s}=P_{A\cap B}y_0.$
\end{remark}

\subsection{Rank-adaptive implementation} \label{subsec:rank_adaptive_imeplementation}

\begin{algorithm}[!b]
\caption{Rank-adaptive $\LSM$ with restarts ($\RaDLSM$)}
\label{alg:lt_rank_detect}
\begin{algorithmic}[1]
\Require Initial point $y_0\in\mathbb R^n$, rank tolerance
$\tau>0$, stopping tolerance $\varepsilon>0$, and maximum
number of DR iterations $N_{\max}$.

\State Set $\nu_0:=0$.
\For{$k=1,2,\ldots$}
    \State Compute $y_k:=T_{A,B}y_{k-1}$
    and $d_k:=y_k-y_{k-1}$.
    \If{the stopping criterion is satisfied at $y_k$}
        \State \Return $y_k$.
    \EndIf
    \State Set $\nu_k:=\rank_\tau(D_k)$.
    \If{$\nu_k=\nu_{k-1}$}
        \State Set $\tilde s:=\nu_k$ and terminate
        the rank-detection phase.
    \EndIf
    \If{$k=N_{\max}$}
        \State Terminate without convergence.
    \EndIf
\EndFor

\State Set $z:=\LSM_{\tilde s}(y_0)$ using the stored
iterates $y_1,\ldots,y_{\tilde s}$.

\While{fewer than $N_{\max}$ DR iterations have been performed}
    \If{the stopping criterion is satisfied at $z$}
        \State \Return $z$.
    \EndIf
    \State Set $z:=\LSM_{\tilde s}(z)$.
\EndWhile
\State Terminate without convergence.
\end{algorithmic}
\end{algorithm}
The finite convergence result in Theorem~\ref{thm:finite_LT} requires the extrapolation order to equal
\[
    s=\dim\mathcal S=
    \dim\spa\{d_1,d_2,\dots\}.
\]
By Lemma~\ref{lem:difference_space_properties}\textup{(iv)}, the first $s$ difference vectors are linearly independent, while every subsequent difference belongs to their span. Consequently, $\rank(D_k)=\min\{k,s\}$, and the dimension $s$ is identified when the rank first fails to increase.

Since $s$ is not known a priori, we estimate it numerically from the DR difference matrices $D_k$. Specifically, as successive DR differences are generated, we monitor $\nu_k:=\rank_{\tau}(D_k)$, where $\tau>0$ is a prescribed rank tolerance. When the numerical rank first fails to increase, we fix the detected dimension to $\tilde s=\nu_k$.

In finite precision, near-linear dependence among the columns of $D_k$ can affect the detected rank. Moreover, even when the correct dimension is detected, the corresponding $\LSM$ system $D_{\tilde s}^\top W_{\tilde s}\alpha =\dia\ \left(D_{\tilde s}^\top W_{\tilde s}\right)$ may be ill-conditioned. When $\tilde s$ is small and this system is sufficiently well conditioned, a single $\LSM_{\tilde s}$ extrapolation may recover a fixed point. However, as the number of active principal-angle directions increases, the system can become increasingly ill-conditioned, and the finite-convergence behaviour predicted by Theorem~\ref{thm:finite_LT} may not be observed numerically.

This motivates the practical rank-adaptive method $(\RaDLSM)$ with restarts summarised in Algorithm~\ref{alg:lt_rank_detect}. After detecting $\tilde s$, we compute an initial $\LSM_{\tilde s}$ extrapolate using the stored DR iterates. If this point does not satisfy the stopping criterion, we generate a new block of $\tilde s$ DR iterates, starting from the extrapolated point, and compute another $\LSM_{\tilde s}$ extrapolate. In this regime, the extrapolations can substantially reduce the fixed-point residual even when they do not produce finite convergence. The process is repeated using the same detected value of $\tilde s$ until the stopping criterion is satisfied or the iteration budget is exhausted.

Each evaluation of $\LSM_{\tilde s}(z)$ generates $\tilde s$ DR iterates from $z$ and constructs the corresponding extrapolate. For the first evaluation, the iterates generated during the rank-detection phase are reused. To mitigate the ill-conditioning noted above, we express the LSM extrapolate in an orthonormal basis, obtained using a QR factorisation \cite[Section~5.2]{golub2013matrix}. This changes the coordinates used to compute the extrapolate while preserving its defining orthogonality conditions.

Specifically, let $W_p=Q_pR_p$ be a reduced QR factorisation. Substituting $\beta=R_p\alpha$ into the LSM coefficient system \eqref{eqn:LSM_coefficient_system} gives
\begin{equation}\label{eq:QR}
    D_p^\top Q_p\beta=\dia\left(D_p^\top W_p\right),\qquad y^{\LSM_p}=y_0+Q_p\beta.
\end{equation}
The right-hand side remains unchanged because only the representation of the extrapolate has changed. Since $W_p=D_p\Gamma_p$, the original coefficient matrix satisfies $D_p^\top W_p=(D_p^\top D_p)\Gamma_p$, where the Gram factor $D_p^\top D_p$ squares the condition number of $D_p$. In contrast, the columns of $Q_p$ form an orthonormal basis for $\range(W_p)=\range(D_p)$, and hence $\kappa(D_p^\top Q_p)=\kappa(D_p)$, where $\kappa(M):=\sigma_{\max}(M)/\sigma_{\min}(M)$ denotes the spectral condition number, with $\kappa(M)=\infty$ when $M$ is singular. Thus, the QR formulation avoids the additional conditioning associated with the direct coefficient matrix, although it cannot remove ill-conditioning already present in $D_p$. All LSM extrapolates in the numerical experiments are computed using this QR formulation; we solve the first equation in \eqref{eq:QR} for $\beta$ and then form $y^{{\LSM}_p }=y_0+Q_p\beta$.

\section{Related extrapolation and Krylov methods}\label{sec:Krylov}
We now place the $\LSM$ recovery result in the context of established extrapolation and Krylov methods. We describe LeSE and DR-GMRES, identify their common search space with $\LSM$, and use Theorem~\ref{thm:dim_reduction} to express their recovery bounds in terms of the active principal angles.

\subsection{Least-squares extrapolation for DR}
\label{subsec:LeSE}
We consider the least-squares extrapolation step underlying A$^3$DMM \cite{poon2019trajectory,poon2020geometry}, applied to the DR difference sequence. We refer to this scheme as $\LeSE$. It fits a linear recurrence to the observed differences and sums the predicted future differences to obtain an extrapolated point.

Fix an extrapolation order $q\ge1$ and $k\ge q+1$. Define 
\[
    \Delta_{k-1}^{(q)}:=[d_{k-1},d_{k-2},\dots,d_{k-q}]\in \R^{n\times q},\qquad c_k\in \argmin_{c\in \R^q}\|\Delta_{k-1}^{(q)}c-d_k\|.
\]
Let
\[
    C_k:=[c_k,e_1,\dots,e_{q-1}],
    \qquad
    \Delta_k^{(q)}:=[d_k,d_{k-1},\dots,d_{k-q+1}],
\]
where $e_i$ denotes the $i$th standard basis vector of $\R^q$ and $C_k=[c_k]$ when $q=1$. Then
\[
    \Delta_{k-1}^{(q)}C_k
    =
    [\Delta_{k-1}^{(q)}c_k,d_{k-1},\dots,d_{k-q+1}]
    \approx \Delta_k^{(q)},
\]
so that $\Delta_{k}^{(q)}C_k^je_1$ approximates $d_{k+j}$, wherefore the extrapolate $y_k + \Delta_{k}^{(q)}\sum_{j=1}^{\tilde{j}}C_k^je_1$ approximates $y_{k+\tilde{j}}$. When $\rho(C_k)<1$, we can use the fact that the Neumann matrix geometric series converges with $\sum_{j=1}^{\infty}C_k^j = (I-C_k)^{-1}-I$ to extrapolate for infinitely many steps as follows:
\[
    y_k^{\LeSE}
    :=
    y_k+\Delta_k^{(q)}
    \bigl((I-C_k)^{-1}-I\bigr)e_1.
\]
The derivation of this infinite-step extrapolation and its connection with minimal polynomial extrapolation (MPE), up to a shift in the iterates used, are discussed in \cite[Section~4.2]{poon2019trajectory}. For affine linear sequences, MPE recovers the limit exactly when its order equals the degree of the minimal polynomial relative to the initial difference \cite{cabay1976polynomial}. In the present setting, Theorem~\ref{thm:dim_reduction} identifies this degree as $s=2|\Theta_{\mathrm{act}}|+\delta_N$. The following corollary verifies exact recovery for the $\LeSE$ construction above.

\begin{corollary}[Exact recovery by $\LeSE$] \label{cor:finite_lsdr}
    Under the assumptions and notation of Theorem~\ref{thm:finite_LT}, suppose that $s\ge 1$. The $\LeSE$ construction with $q=s$ and $k=s+1$ satisfies
    \[
        \rho(C_{s+1})<1,\qquad y_{s+1}^{\LeSE}=P_{\Fix(T_{A,B})}y_0.
    \]
\end{corollary}

\begin{proof}
By Lemma~\ref{lem:difference_space_properties}, the columns of $\Delta_s^{(s)}$ form a basis of $\mathcal S$. Since $d_{s+1}\in \mathcal S$, the least-squares fit is exact, and the difference recurrence gives $L\Delta_s^{(s)}=\Delta_s^{(s)}C_{s+1}$. Thus, $C_{s+1}$ represents $L|_S$. Since $\mathcal S\subseteq G$ and every eigenvalue of $L|_G$ has modulus less than one, $\rho(C_{s+1})<1$.

Moreover, for every $j\ge1$, 
\[
    \Delta_{s+1}^{(s)}C_{s+1}^{j}e_1=d_{s+1+j}.
\]
Expanding $(I-C_{s+1})^{-1}-I$ as a convergent geometric series yields
\begin{align*}
    y_{s+1}^\LeSE
    &=y_{s+1}+\Delta_{s+1}^{(s)}((I-C_{s+1})^{-1}-I)e_1\\
    &=y_{s+1}+ \Delta_{s+1}^{(s)}\sum_{j=1}^{\infty}C_{s+1}^je_1\\
    &=y_{s+1}+\sum_{j=1}^\infty d_{s+1+j}\\
    &=y_{s+1}+\lim_{N\to\infty}\sum_{j=1}^N(y_{s+1+j}-y_{s+j})\\
    &=y_{s+1}+\lim_{N\to\infty}(y_{s+1+N}-y_{s+1})=\lim_{k\to\infty}y_k=P_{\Fix(T_{A,B})}y_0,
\end{align*}
where the final equality follows from Proposition~\ref{prop:affine_DR_limit}.
\end{proof}

Given a starting point $z\in \R^n$, let $\LeSE_q(z)$ denote the point obtained by generating $q+1$ DR differences from $z$ and applying the preceding construction with $k=q+1$. Corollary~\ref{cor:finite_lsdr} guarantees $\rho(C_{q+1})<1$ when $q$ equals the difference-space dimension. In implementation, however, near-linear dependence among the observed differences can lead to inaccurate recurrence coefficients, and the computed matrix may not satisfy this condition. We therefore check it explicitly and retain the most recent DR iterate if it fails.

Corollary~\ref{cor:finite_lsdr} motivates using the detected dimension $\tilde s$ as the $\LeSE$ extrapolation order. The resulting rank-adaptive method, $\RaDLeSE$, follows the same dimension-detection and restart framework as $\RaDLSM$. Once $\tilde s$ is detected, we set $q=\tilde s$ and reuse the stored differences to compute the first extrapolate. If the stopping criterion is not satisfied, we apply $\LeSE_{\tilde s}$ repeatedly, starting each block from the point returned by the preceding block. The detected order remains fixed until the stopping criterion is satisfied or the iteration budget is exhausted. Each restarted $\LeSE$ block requires $\tilde s+1$ DR differences, compared with $\tilde s$ for an $\LSM _{\tilde s}$ block. The additional difference is needed to fit the LeSE recurrence coefficients using the preceding $\tilde s$ differences.

\subsection{GMRES applied to the DR fixed-point equation} \label{subsec:GMRES}
The generalized minimal residual method (GMRES) \cite{saad1986gmres} is a Krylov subspace method for solving linear systems. For a system $Mx=b$ and an initial point $x_0$, let $r_0:=b-Mx_0$ be the initial residual, and define the associated Krylov subspace of order $k$ as $\mathcal K_k(M,r_0):=\spa\{r_0,Mr_0,\dots,M^{k-1}r_0\}$. At iteration $k$, GMRES chooses
\[
x_k\in \argmin_{x\in x_0+\mathcal K_k(M,r_0)}\|b-Mx\|.
\]
The method constructs an orthonormal basis of the Krylov space using the Arnoldi process and computes the approximation by solving a reduced least-squares problem.

For affine DR, $T_{A,B}y=Ly+c$, so the fixed-point equation is equivalent to $(I-L)y=c$. Applying GMRES to this system from $y_0$ gives the initial residual $r_0=c-(I-L)y_0=T_{A,B}y_0-y_0=d_1$, and the approximations satisfy
\[
    y_k^\gmres\in\argmin_{y\in y_0+\mathcal K_k(I-L,r_0)}\|T_{A,B}y-y\|.
\]
Although $I-L$ may be singular on $\R^n$, the initial residual belongs to $G$, which is invariant under both $L$ and $I-L$. Consequently, the Krylov subspaces lie in $G$, where $(I-L)|_G$ is nonsingular.

We apply GMRES to the DR fixed-point equation to investigate its relationship with extrapolation based on successive DR iterates. This formulation allows us to compare the search spaces and numerical behaviour of $\gmres$, $\LSM$, and $\LeSE$.

\subsection{A common search space} \label{subsec:common_search_space}
The Krylov space used by $\gmres$ coincides with the space spanned by the observed DR differences. The following proposition makes this relationship explicit.

\begin{proposition}[Common search space] \label{prop:common_search_space}
Let $T_{A,B}y=Ly+c$, let $y_{j+1}=T_{A,B}y_j$, and set $r_0=T_{A,B}y_0-y_0=d_1$. Then, for every $k\ge1$,
\[
    \mathcal K_k(I-L,r_0)
    =\mathcal K_k(L,d_1)
    =\spa\{d_1,\dots,d_k\}
    =\range(D_k)
    =\range(W_k).
\]
Consequently, $y_0+\mathcal K_k(I-L,r_0)=\aff\{y_0,\dots,y_k\}$.
\end{proposition}

\begin{proof}
Since $r_0=d_1$, the binomial expansion of $(I-L)^j$ gives $(I-L)^jd_1\in\spa\{d_1,Ld_1,\dots, L^jd_1\}$. Hence $\mathcal K_k(I-L,r_0)\subseteq\mathcal K_k(L,d_1)$. Conversely, expanding $L^j=(I-(I-L))^j$ gives the reverse inclusion.

The recurrence $d_j=L^{j-1}d_1$ yields $\mathcal K_k(L,d_1)=\spa\{d_1,\dots,d_k\}=\range(D_k)$. Moreover, $W_k=D_k\Gamma_k$, where $\Gamma_k$ is the nonsingular matrix defined in Section~\ref{subsec:LSM_extrap}. Therefore $\range(W_k)=\range(D_k)$. Finally, $\aff\{y_0,\dots,y_k\}=y_0+\range(W_k)$.
\end{proof}

When constructed from $y_0,\dots,y_k$, the $\LSM$ extrapolate, whenever defined, belongs to the affine hull of the iterates used to construct it. The $\LeSE$ extrapolate belongs to the same affine hull, since it is obtained by adding a linear combination of the observed differences to $y_k$. The same holds when the most recent DR iterate $y_k$ is retained instead. By Proposition~\ref{prop:common_search_space}, this affine hull is the search space used by $\gmres$ at iteration $k$. Thus, all three methods select points from a common affine space using different criteria.

Finite termination of GMRES for nonsingular linear systems is classical \cite{saad1986gmres}. Here, the common-space identification shows that $\gmres$ recovers the projected DR fixed point within $s$ iterations.

\begin{corollary}[Exact recovery by $\gmres$] \label{cor:finite_GMRES}
Under the assumptions and notation of Theorem~\ref{thm:finite_LT}, apply $\gmres$ to $(I-L)y=c$ with initial point $y_0$. The method terminates in at most $s$ iterations and returns $P_{\Fix(T_{A,B})}y_0$.    
\end{corollary}

\begin{proof}
    Set $\bar y:=P_{\Fix(T_{A,B})}y_0$. If $s=0$, then $d_1=0$, so $y_0=\bar y$. Suppose that $s\ge 1$. By \eqref{eq:projection_in_difference_space} and Proposition~\ref{prop:common_search_space}, $\bar y-y_0\in \mathcal S=\mathcal K_s(I-L,r_0)$. Thus, $\bar y$ belongs to the $\gmres$ search space at iteration $s$ and has zero residual. Consequently, $\gmres$ terminates in at most $s$ iterations.

    The returned point belongs to $y_0+\mathcal S$ and has zero residual. Since $\mathcal S\subseteq G$ and $(I-L)|_G$ is nonsingular, $\bar y$ is the unique solution of $(I-L)y=c$ in $y_0+\mathcal S$. Therefore, the returned point is $\bar y$.
\end{proof}

Theorem~\ref{thm:finite_LT} and Corollaries~\ref{cor:finite_lsdr} and \ref{cor:finite_GMRES} show that full-order $\LSM$, $\LeSE$, and $\gmres$ recover the same fixed point $P_{\Fix(T_{A,B})}y_0$. Their numerical behaviour, however, depends on how the common search space is represented.

For $k\le s$, the columns of $D_k=[d_1,Ld_1,\dots, L^{k-1}d_1]$ form a power basis for $\mathcal K_k(L,d_1)$. These vectors may become nearly linearly dependent, making the basis poorly conditioned \cite[Section~7.3.1]{hoemmen2010communication}. $\LSM$ and $\LeSE$ compute their coefficients from the observed differences. The QR formulation described in Section~\ref{subsec:rank_adaptive_imeplementation} improves the coordinates used in the $\LSM$ solve, but does not change the difference matrix used for dimension detection. $\gmres$ instead constructs an orthonormal basis incrementally using the Arnoldi process \cite{saad1986gmres}. Thus, the methods can behave differently in finite precision even though their search spaces coincide in exact arithmetic.
\color{black}
\section{Numerical experiments} \label{sec:numerical_experiments}
The numerical experiments have three objectives. First, we verify the difference-space dimension formula $s=2|\Theta_{\mathrm{act}}|+\delta_N$ from Theorem~\ref{thm:dim_reduction} and the full-order $\LSM$ recovery result of Theorem~\ref{thm:finite_LT}. Second, we examine the ill-conditioning effects discussed in Section~\ref{subsec:rank_adaptive_imeplementation} when principal angles are nearly equal. Third, we compare the practical performance of $\RaDLSM$ and $\RaDLeSE$ with DR, $\gmres$, and fixed-order variants of $\LSM$ and $\LeSE$ as the number of active principal angles increases.
\subsection{Experimental setup}
We consider consistent affine subspaces
\[
A=x_\star+U\qquad \text{and}\qquad B=x_\star+V
\]
of $\R^{400}$. Unless otherwise stated, their direction subspaces satisfy
\[
\dim(U\cap V)=50,\qquad \dim U=120,\qquad \dim V=160.
\]
The subspaces are generated according to the principal-angle structure of Lemma~\ref{lem:principal_structure}. We prescribe the distinct proper principal angles and the dimensions of the components in the orthogonal decomposition given there. The resulting decomposition is expressed in a random orthonormal basis, and the common affine point $x_\star$ is sampled from a standard Gaussian distribution.

Numerical ranks are computed using the relative singular-value threshold
\[
    \rank_\tau(D_k):=\bigl|\{i:\sigma_i(D_k)>\tau\sigma_1(D_k)\}\bigr|,
\]
where $\sigma_1(D_k)\ge \sigma_2(D_k)\ge\cdots$ denote the singular values of $D_k$. Unless otherwise stated, we take $\tau=10^{-13}$. The detected dimension $\tilde s$ is fixed when the numerical rank fails to increase. 

We use the fixed point-residual $r_{\Fix}(y)=\|T_{A,B}y-y\|.$ A run is declared successful when $r_\Fix(y)\le10^{-6}$ within a maximum of $10,000$ operator evaluations. Runtimes include rank detection, extrapolation calculations, and stopping criterion evaluations, but exclude problem generation and plotting.

Reported operator counts include evaluations used to generate the iterates and extrapolation data. For $\gmres$, we count the initial residual evaluation and each application of $I-L$. 
\subsection{Difference-space dimension and conditioning}
We examine the conditioning of the direct LSM coefficient matrix $D_p^\top W_p$ and its QR counterpart $D_p^\top Q_p$, introduced in Section~\ref{subsec:rank_adaptive_imeplementation}, and investigate how near-linear dependence among the DR differences affects dimension detection and exact recovery.

\paragraph{Controlled active blocks.}
For this experiment, we prescribe four distinct principal angles $\Theta=\{10^\circ,25^\circ,40^\circ,55^\circ\}$ and a nontrivial null component $N$. Proposition~\ref{prop:deg_T} therefore gives the full difference-space bound $\ell_G=2|\Theta|+1=9$. For each $j\in \{0,\dots,4\}$ and $\delta_N\in \{0,1\}$, we select a subset $I\subseteq\Theta$ with $|I|=j$ and choose an initial displacement $g\in G$ satisfying
\[
    P_{W_\theta}g\neq0\qquad \longleftrightarrow\qquad \theta\in I,\qquad P_Ng\neq0\quad \longleftrightarrow\quad \delta_N=1,
\]
with all remaining components equal to zero. The nonzero components are generated randomly, and $g$ is normalised whenever $g\ne 0$. We then set $y_0=x_\star+g.$ By construction, $\Theta_{\mathrm{act}}=I$, and the null component is active precisely when $\delta_N=1$. Theorem~\ref{thm:dim_reduction} therefore predicts $s=2|I|+\delta_N$. Moreover, because $g\in G=F^\perp$ and $\Fix(T_{A,B})=x_\star+F$, the distinguished fixed point is known explicitly:
\[
    P_{\Fix(T_{A,B})}y_0=x_\star.
\]

For each pair $(|I|,\delta_N)$, we generate 25 randomly oriented instances and compute the first ten DR differences. The detected dimension $\tilde s$ is obtained from the first plateau in the numerical ranks of the growing difference matrices. For $s>0$, we also compute the full-order $\LSM_s$ extrapolated point using the value predicted by Theorem~\ref{thm:dim_reduction} and assess the recovery of Theorem~\ref{thm:finite_LT} using $e_\star(z_s^\LSM):=\|z_s^\LSM-x_\star\|$, together with the fixed-point residual $r_{\Fix}(z_s^\LSM)$.

\begin{table}
\centering
\caption{Controlled activation of the principal-angle and null components. The residual and target-error columns report medians over 25 trials.}
\label{tab:dimension_bound}
\begin{tabular}{cccccc}
\toprule
$|\Theta_{\mathrm{act}}|$ & $\delta_N$ & Predicted $s$ &
Detected $\tilde{s}$ & Median $r_{\mathrm{Fix}}$ &
Median $e_\star$ \\
\midrule
0 & 0 & 0 & 0 & $0$ & $0$ \\
0 & 1 & 1 & 1 & $7.25\times 10^{-16}$ & $1.06\times 10^{-15}$ \\
1 & 0 & 2 & 2 & $1.07\times 10^{-14}$ & $1.28\times 10^{-14}$ \\
1 & 1 & 3 & 3 & $4.29\times 10^{-15}$ & $5.24\times 10^{-15}$ \\
2 & 0 & 4 & 4 & $5.68\times 10^{-14}$ & $6.90\times 10^{-14}$ \\
2 & 1 & 5 & 5 & $1.64\times 10^{-13}$ & $1.92\times 10^{-13}$ \\
3 & 0 & 6 & 6 & $3.21\times 10^{-12}$ & $3.68\times 10^{-12}$ \\
3 & 1 & 7 & 7 & $3.01\times 10^{-12}$ & $3.57\times 10^{-12}$ \\
4 & 0 & 8 & 8 & $1.93\times 10^{-11}$ & $2.25\times 10^{-11}$ \\
4 & 1 & 9 & 9 & $1.97\times 10^{-11}$ & $2.29\times 10^{-11}$ \\
\bottomrule
\end{tabular}
\end{table}
Table~\ref{tab:dimension_bound} shows that the predicted and detected dimensions agree for every activation pattern in all 25 trials. Each active principal angle contributes two dimensions, whereas activation of the null component contributes one. In particular, the full bound $\ell_G=9$ is attained precisely when all four principal angle blocks and the null component are active, in agreement with Corollary~\ref{cor:dimension_bound_attainment}.

The full-order $\LSM$ extrapolates also recover the distinguished fixed point to high numerical accuracy. In the full-dimensional case $s=9$, the median fixed-point residual and target error are $1.97\times10^{-11}$ and $2.29\times10^{-11}$, respectively. The largest errors occur at the higher extrapolation orders, although they remain small relative to the stopping tolerance.

\paragraph{Nearly coincident principal angles.}
We keep four principal-angle blocks and the null component active while decreasing the separation between consecutive angles. Specifically, we set $\theta_j(\delta)=30^\circ+(j-\frac{5}{2})\delta$, for $j=1,\dots,4$, where $\delta$ is measured in degrees. We consider $\delta\in\{10,5,2,1,0.5,0.2,0.1,0.05,0.02,0.01,0\}$. For every $\delta>0$, the four angles are distinct. Since all four principal-angle blocks and the null component are active, Theorem~\ref{thm:dim_reduction} gives $s=9$. At $\delta=0$, the four angles coincide, so there is only one distinct active principal angle and the exact dimension becomes $s=3$.

For each separation, we generate 25 instances. Within each trial, the principal planes, the remaining subspace components, the common point $x_\star$, and the initial displacement $g=y_0-x_\star$ are held fixed as $\delta$ varies. Only the angles within the principal planes are changed. For each $\delta>0$, we compute the condition numbers of the direct and QR-reparameterised $\LSM$ systems, $\kappa(D_9^\top W_9 )\text{ and } \kappa(D_9^\top Q_9),$ respectively. We also apply the rank detector using $\tau\in\{10^{-8},10^{-10},10^{-12},10^{-13}\}$. All reported quantities are medians over the 25 trials.

\begin{figure}[h]
    \centering
    \includegraphics[width=\linewidth]{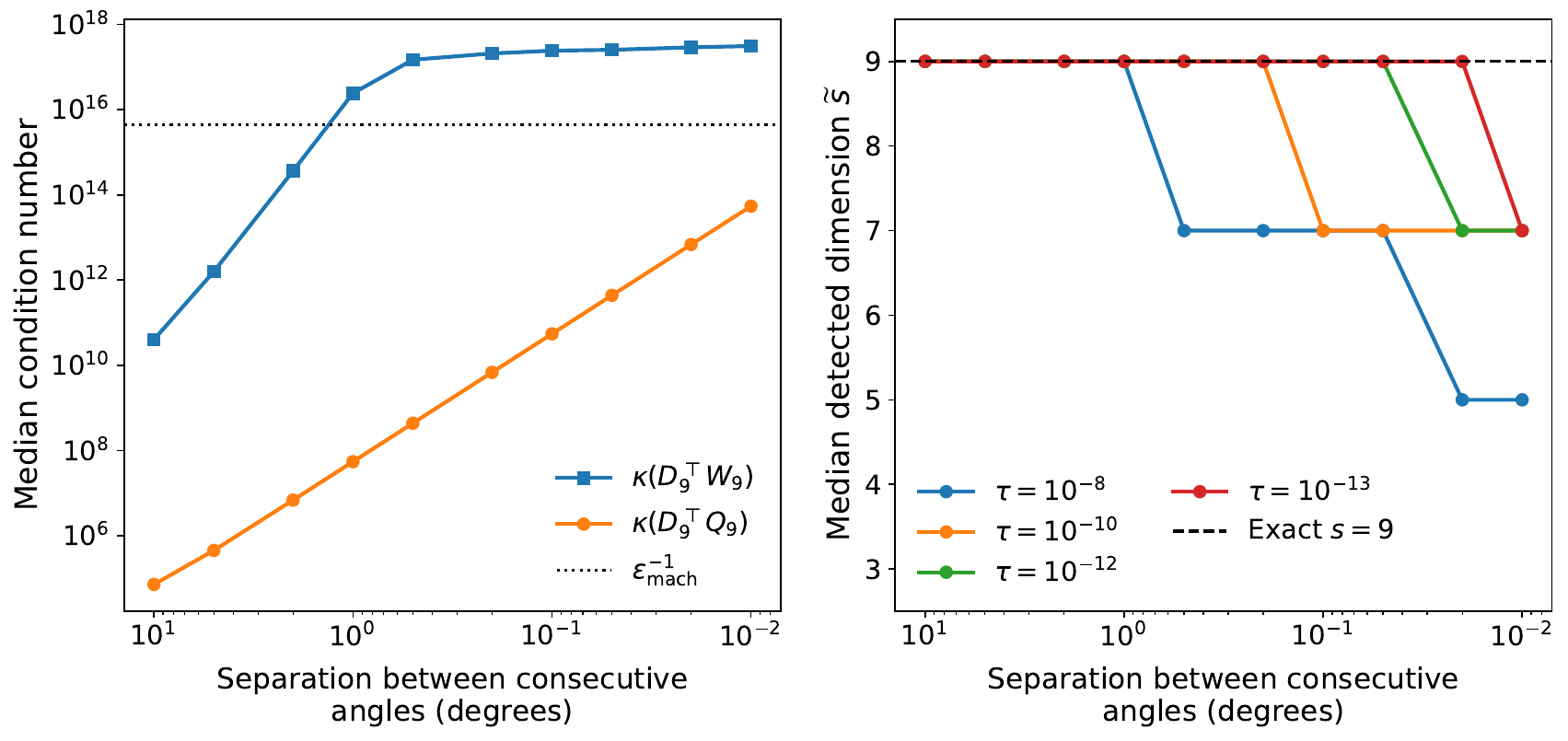}
    \caption{Effect of decreasing principal-angle separation on conditioning and dimension detection. Left: median condition numbers of the direct and QR-reparameterised $\LSM$ systems. Right: median detected dimension for several rank tolerances; the dashed line marks the exact dimension $s=9$.}
    \label{fig:conditioning}
\end{figure}
Figure~\ref{fig:conditioning} (left) shows that both coefficient matrices become increasingly ill-conditioned as the separation decreases, but at substantially different rates. The direct matrix $D_9^\top W_9$ reaches the scale of inverse machine precision much earlier than the QR-reparameterised matrix $D_9^\top Q_9$. The QR formulation therefore removes much of the additional ill-conditioning introduced by the nonorthogonal coordinates in $W_9$, although it cannot remove the ill-conditioning already present in $D_9$ as the angles approach one another.

Figure~\ref{fig:conditioning} (right) shows how the detected dimension depends on the rank tolerance. With $\tau=10^{-8}$, the median detected dimension decreases from nine to seven at $\delta=0.5^\circ$ and then to five at $\delta=0.02^\circ$. With $\tau=10^{-13}$, it remains equal to nine through $\delta=0.02^\circ$ and decreases to seven only at $\delta=0.01^\circ$. Thus, smaller rank tolerances distinguish the nearly coincident components over a wider range of separations. The exact dimension remains $s=9$ throughout the figure because the four angles are distinct at every plotted separation.

\begin{table}
\centering
\caption{Selected separations for four active principal-angle components. Numerical ranks use $\tau=10^{-13}$.}
\label{tab:nearly_equal_angles}
\begin{tabular}{cccccc}
\toprule
Separation $(^\circ)$ & $|\Theta_{\mathrm{act}}|$ & Exact $s$ & Detected $\tilde{s}$ & $\kappa(D_s^\top W_s)$ & $\kappa(D_s^\top Q_s)$ \\
\midrule
0.02 & 4 & 9 & 9 & $2.90\times 10^{17}$ & $6.86\times 10^{12}$ \\
0.01 & 4 & 9 & 7 & $3.11\times 10^{17}$ & $5.44\times 10^{13}$ \\
0 & 1 & 3 & 3 & $1.74\times 10^{2}$ & $6.64\times 10^{0}$ \\
\bottomrule
\end{tabular}
\end{table}

Table~\ref{tab:nearly_equal_angles} isolates the transition at the default rank tolerance. At $\delta=0.02^\circ$, the detector returns the exact dimension $s=9$. However, the direct $\LSM$ matrix has condition number $2.90\times10^{17}$, whereas the QR-reparameterised matrix has condition number $6.86\times 10^{12}$. The QR formulation therefore improves the conditioning considerably at this separation. At $\delta=0.01^\circ$, the condition number of the QR system increases to $5.44\times10^{13}$. Equivalently, its smallest relative singular value is approximately $1.84\times 10^{-14}$, below the rank threshold $\tau=10^{-13}$. This is consistent with the detector resolving only $\tilde s=7$ directions.

At exact equality, the reduction to $s=3$ is algebraic rather than numerical. The four angle components share the same quadratic factor and collectively contribute two dimensions, while the active null component contributes one. The detector therefore correctly returns $\tilde s=3$, and the corresponding full-order matrices are substantially better conditioned.

\subsection{Practical performance}
\label{subsec:practical_peformance}
We investigate how numerical dimension detection affects the practical performance of the extrapolation methods as the number of distinct principal angles increases. We consider $m\in\{1,2,4,6,8,10,15,20,30,40,50\}$. For $m\ge2$, the angles are equally spaced over $[5^\circ,60^\circ]$; for $m=1$, the single angle is $5^\circ$. The principal-angle decomposition is expressed in a random orthonormal basis, and the initial points are sampled from a standard Gaussian distribution. Consequently, all principal-angle components and the null component are active, and Theorem~\ref{thm:dim_reduction} predicts $s=2m+1$.

For each $m$, we generate $25$ independent affine-subspace pairs and use $5$ initial points per pair. We compare DR, $\LSM_2$, $\LeSE_6$, $\RaDLSM$, $\RaDLeSE$, and $\gmres$ on the same $125$ problems. We report median DR operator counts and runtimes.

\paragraph{Detected and Arnoldi dimensions.}
\begin{figure}[t]
    \centering 
    \includegraphics[width=0.72\linewidth]{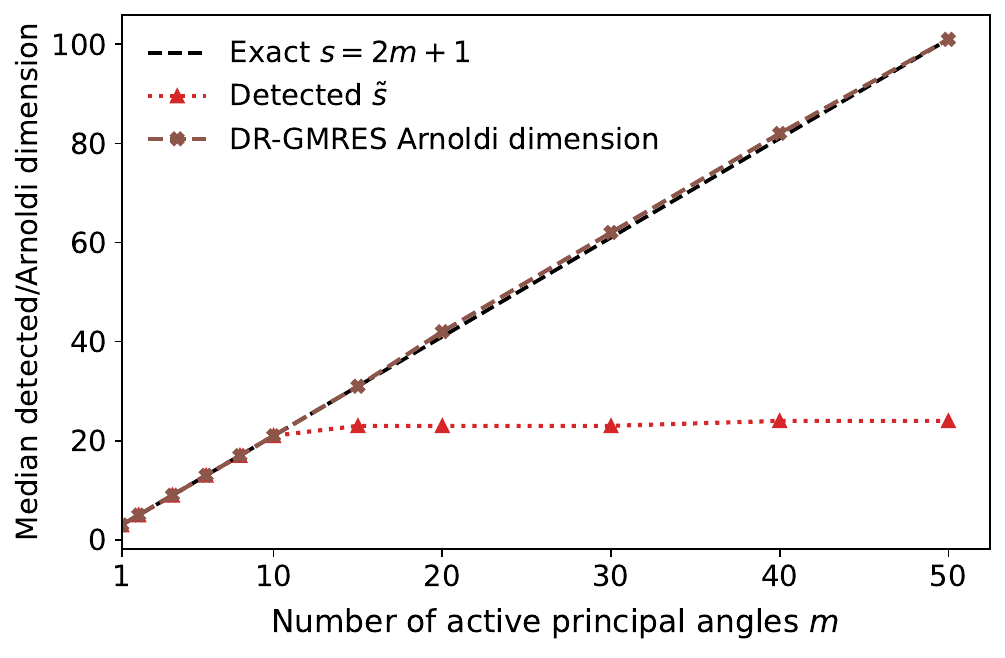} 
    \caption{Exact difference-space dimension, median rank-detected extrapolation
    order, and median $\gmres$ Arnoldi dimension as the number \(m\) of active
    principal angles increases.}
    \label{fig:working_dimension} 
\end{figure}
$\RaDLSM$ and $\RaDLeSE$ apply the same rank detector to the same initial DR differences, so they select the same order $\tilde s$. Figure~\ref{fig:working_dimension} compares this detected order with the exact dimension $s$ and the number of Arnoldi directions generated by $\gmres$ before reaching the stopping tolerance.

The median detected order equals $s$ through $m=10$. At $m=15$, it falls to $\tilde s=23$, compared with $s=31$, and remains approximately $23$ or $24$ as $m$ increases further. This behaviour is consistent with the conditioning effects observed in the preceding subsection: placing more angles in the fixed interval $[5^\circ,60^\circ]$ reduces their separation, making independent DR differences harder to distinguish numerically.

In contrast, $\gmres$ reaches the stopping tolerance in a single Arnoldi process on every problem. Its median Arnoldi dimension equals $s$ for $m\le15$, and $m=50$, and is $s+1$ for $m=20,\,30,\, 40$. Corollary~\ref{cor:finite_GMRES} guarantees termination by step $s$, but rounding errors can require additional Arnoldi directions to reach the prescribed tolerance.

\paragraph{Increasing difference-space dimension.}
\begin{figure}[t]  
    \centering
    \includegraphics[width=\linewidth]
    {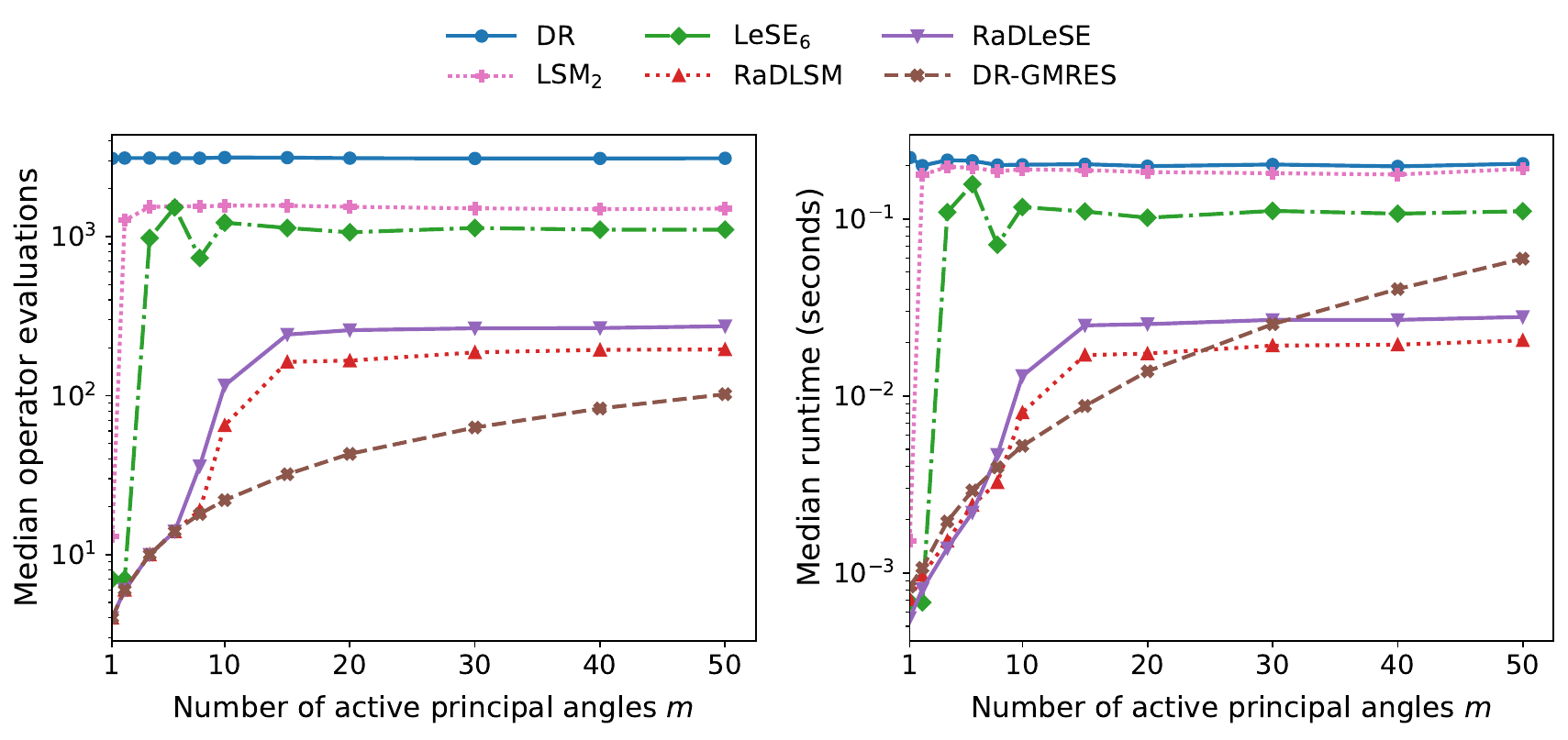}
    \caption{Median operator evaluations and runtime as the number \(m\) of
    active principal angles increases. Each point is the median over
    25 affine-subspace pairs and five initial points per pair.}
    \label{fig:cost_vs_m}
\end{figure}

Figure~\ref{fig:cost_vs_m} (left) shows that both rank-adaptive methods have the same median DR operator count as $\gmres$ for $m\le6$. At $m=8$ and $m=10$, they require additional DR steps or further extrapolations despite detecting the correct order. Thus, correct dimension detection alone does not ensure accurate recovery in a single extrapolation.

For larger values of $m$, the adaptive methods apply extrapolation repeatedly at the smaller detected order, whereas $\gmres$ continues expanding its Krylov space. The adaptive methods require more DR evaluations, but their counts begin to level off as $\tilde s$ stabilises. At $m=50$, the median counts are $102$ for $\gmres$, $195$ for $\RaDLSM$, and $273$ for $\RaDLeSE$. Although the adaptive methods use the same order, they can require different numbers of extrapolations. Each new $\LSM_{\tilde s}$ extrapolate requires $\tilde s$ steps, compared with $\tilde s+1$ for $\LeSE_{\tilde s}$. 

The fixed-order methods improve on unmodified DR but require approx. $10^3$ DR evaluations over much of the tested range. Their prescribed orders become insufficient for exact recovery as $s$ increases, so repeated extrapolation is needed. Unmodified DR requires approximately $3.1\times10^3$ evaluations throughout. Its nearly constant count is consistent with the fixed Friedrichs angle of $5^\circ$, which determines its asymptotic convergence rate \cite{BAUSCHKE201463}.

Figure~\ref{fig:cost_vs_m} (right) shows that the smallest DR operator count does not always give the shortest runtime. $\gmres$ has the lowest median runtime for $m=10,15,20$. However, $\RaDLSM$ is faster for $m\ge30$, and $\RaDLeSE$ is faster for $m\ge 40$. At $m=50$, their median runtimes are approximately $0.017$ and $0.024$ seconds, respectively, compared with $0.053$ seconds for $\gmres$.

This change is consistent with the different computational costs of the methods. The adaptive methods retain a small detected order and repeatedly solve systems of that size, while $\gmres$ constructs and orthogonalises a growing Arnoldi basis. The fixed-order methods also solve small systems, but require sufficiently many repeated extrapolations that they do not achieve comparable runtimes over most of the tested range.

\section{Conclusion}\label{sec:conclusion}
For two intersecting affine subspaces, the dimension of the affine space generated by the DR iterates depends on the distinct nontrivial principal angles active in the initial difference. A low extrapolation order can therefore suffice for exact recovery, even in high-dimensional problems. We have shown that LSM, LeSE, and DR-GMRES use the same affine search spaces and recover the same DR limit from finitely many iterates, with recovery bounds determined by the dimension of the space.

Our rank-adaptive implementations of LSM and LeSE estimate the extrapolation order from the observed DR differences. Nearly equal principal angles can produce almost linearly dependent DR differences, leading to an underestimate of the required extrapolation order. Even when the correct order is detected, the extrapolation system can be ill-conditioned. The QR formulation improves the numerical stability of the LSM solve, although ill-conditioning in the difference matrix can still limit its accuracy.

The rank-adaptive methods still required fewer DR operator evaluations than unmodified DR when the extrapolation order was underestimated. DR-GMRES required the fewest evaluations in the tested problems, but the adaptive methods sometimes ran faster at larger search-space dimensions. This difference in runtime is consistent with the adaptive methods repeatedly solving smaller systems, while DR-GMRES constructs a growing Arnoldi basis.

The affine analysis also provides a starting point for studying extrapolation in nonlinear feasibility problems. In applications such as phase retrieval, smooth constraint sets can be approximated near a solution by their tangent affine spaces. A natural question is whether the geometric structure underlying finite recovery in the affine setting can be used to establish local convergence rates for extrapolation methods on curved constraint sets.

\bibliographystyle{plainnocase}
\bibliography{bibliography}

\appendix
\section{Proof of Proposition~\ref{prop:angles}}
\label{app:dr_proof}

\begin{proof}
    Fix a principal angle $\theta_i\in(0,\tfrac{\pi}{2})$. By Lemma~\ref{lem:princ_planes}, $(u_i,\widetilde v_i)$ is an orthonormal basis of $\Pi_i$, with $v_i=\cos\theta_i\,u_i+\sin\theta_i\,\widetilde v_i$. Since $\widetilde v_i\in U^\perp$, every $x\in\Pi_i$ satisfies $P_Ux\in \spa \{u_i\}\subseteq\Pi_i$. Moreover, Lemma~\ref{lem:principal_structure} shows that $\Pi_i$ is orthogonal to every component of $V$ other than $\spa\{v_i\}$. Hence 
    \[
    P_Vx\in\spa\{v_i\}\subseteq \Pi_i\qquad (x\in \Pi_i).
    \] 
    Thus $P_U$ and $P_V$ leave $\Pi_i$ invariant. Consequently, $R_U$, $R_V$, and $L=\tfrac{1}{2}(I+R_VR_U)$ also leave $\Pi_i$ invariant.

    With respect to the orthonormal basis $(u_i,\widetilde v_i)$, we have
    \[
        P_U|_{\Pi_i}
        =
        \begin{bmatrix}
            1&0\\
            0&0
        \end{bmatrix},
        \qquad
        P_V|_{\Pi_i}
        =
        \begin{bmatrix}
            \cos^2\theta_i
            &
            \sin\theta_i\cos\theta_i\\
            \sin\theta_i\cos\theta_i
            &
            \sin^2\theta_i
        \end{bmatrix}.
    \]
    Therefore,
    \[
        R_U|_{\Pi_i}
        =
        \begin{bmatrix}
            1&0\\
            0&-1
        \end{bmatrix},
        \qquad
        R_V|_{\Pi_i}
        =
        \begin{bmatrix}
            \cos(2\theta_i)&\sin(2\theta_i)\\
            \sin(2\theta_i)&-\cos(2\theta_i)
        \end{bmatrix},
    \]
    and consequently
    \[
        (R_VR_U)|_{\Pi_i}
        =
        \begin{bmatrix}
            \cos(2\theta_i)&-\sin(2\theta_i)\\
            \sin(2\theta_i)&\cos(2\theta_i)
        \end{bmatrix}
        =\Rot(2\theta_i).
    \]
    It follows that
    \begin{align*}
        L|_{\Pi_i}
        &=
        \frac12\bigl(I+(R_VR_U)|_{\Pi_i}\bigr)\\
        &=
        \frac12
        \begin{bmatrix}
            1+\cos(2\theta_i)&-\sin(2\theta_i)\\
            \sin(2\theta_i)&1+\cos(2\theta_i)
        \end{bmatrix}\\
        &=
        \cos\theta_i
        \begin{bmatrix}
            \cos\theta_i&-\sin\theta_i\\
            \sin\theta_i&\cos\theta_i
        \end{bmatrix}\\
        &=
        \cos\theta_i\Rot(\theta_i).
    \end{align*}

    It remains to consider the four boundary components. On $U\cap V$, both reflections act as the identity, whereas on $U^\perp\cap V^\perp$ both act as $-I$. Thus $R_VR_U=I$ and hence $L=I$ on both components.

    On $U\cap V^\perp$, the reflections $R_U$ and $R_V$ act as $I$ and $-I$, respectively. On $U^\perp\cap V$, they act as $-I$ and $I$, respectively. Hence $R_VR_U=-I$ and therefore $L=0$ on both components.
\end{proof}
\end{document}